# PALM DISTRIBUTIONS OF WAVE CHARACTERISTICS IN ENCOUNTERING SEAS

By Sofia Aberg,[1] Igor Rychlik and M. Ross Leadbetter

*Lund University, Lund University and University of North Carolina*

Distributions of wave characteristics of ocean waves, such as wave slope, waveheight or wavelength, are an important tool in a variety of oceanographic applications such as safety of ocean structures or in the study of ship stability, as will be the focus in this paper. We derive Palm distributions of several wave characteristics that can be related to steepness of waves for two different cases, namely for waves observed along a line at a fixed time point and for waves encountering a ship sailing on the ocean. The relation between the distributions obtained in the two cases is also given physical interpretation in terms of a "Doppler shift" that is related to the velocity of the ship and the velocities of the individual waves.

**1. Introduction.** The study of wave characteristics such as crest height, wavelength and wave slope is important in various ocean engineering applications. Examples can be found in the design of ocean structures, such as oil platforms, or in the design of sea walls that should prevent cities, or sometimes whole countries, from flooding. Another important application, motivating this paper, is that of ship stability. Because high and steep waves encountering a ship may cause structural damage and even capsize smaller vessels, it is important to know the distribution of, for example, wavelength, waveheight and wave slope of such waves. Based on these distributions, capsize probabilities can be computed and used as a risk measure for existing vessels or as a tool in the design of new ones.

A sea surface can be seen as a sequence of apparent waves. By an apparent wave is meant the part of the sea record between two consecutive upcrossings of the still water level and the downcrossing in between these upcrossings is called the center of the wave. The purpose of this paper is to derive exact

Received May 2006; revised June 2007.

[1]Supported by the Swedish Foundation for Strategic Research Grant A3 02:125.

distributions of steepness related wave characteristics of the apparent waves for two different cases. The first case, *the spatial case*, concerns properties of the waves observed, in space, along a fixed line at a fixed time point and the second case, the *encountered case*, involves properties of the waves encountering a ship on the ocean. More precisely, the characteristics of those waves that overtake a ship sailing with constant velocity along a straight line will be investigated because such waves are considered as particularly dangerous for ship stability.

Distributions of wave characteristics are defined as Palm distributions. Such distributions are intimately related to level-crossings of the underlying process; in our case representing the sea surface, and Rice's formula proves useful in their evaluation. Although the main result of this paper is the derivation of several distributions for the two different cases stated above, a nice by-product is that the relation between the two cases can be given interpretation in terms of wave velocities. In physical terminology this can be expressed as a Doppler shift and is caused by the fact that the sea surface is observed by a moving observer—the ship. Due to this interpretation, similar results can be expected in many other applications where the distribution of moving objects, observed from a moving observer are studied. As an example, the distribution of the size of the storms you will meet when sailing will be related, by a Doppler shift, to the distribution of storm sizes that can be observed from a satellite image.

The organization of this paper is as follows. First, the Gaussian sea model used to evaluate the distributions for the wave characteristics is introduced. Then it is shown how distributions defined by Palm distributions can be computed by using generalized versions of Rice's formula. After a short section on wave velocities, the theory is exemplified by calculation of the distribution of wave slope for spatial and encountered waves. Thereafter, a more intricate example concerning waveheight and wavelength of apparent waves is investigated, and finally, the derived distributions are evaluated numerically for a unidirectional Gaussian sea.

**2. Gaussian sea model.** Let $W(x,t)$ be the sea surface elevation at location $x$ and time $t$. Further denote its partial derivative on $t$ by $W_t(x,t)$ and its first and second partial derivative on $x$ by $W_x(x,t)$ and $W_{xx}(x,t)$, respectively. Sometimes, when no misunderstanding can be made, the notation $W(x) = W(x,0)$ is used. Throughout this paper, $W(x,t)$ is modeled as a zero mean, stationary Gaussian field with a directional spectrum having spectral density $S(\omega, \theta)$, where $\omega$ is an angular frequency and $\theta$ an angle representing wave direction. These quantities satisfy $\omega > 0$ and $\theta \in [0, 2\pi]$. For sea waves such a parametrization of the spectral density is possible because the angular frequency $\omega$ and the wave number $\kappa$ are related by physical dispersion laws; see [4] for further reference. In particular, for deep water,

$$(2.1) \qquad \kappa = \frac{\omega^2}{g},$$

where $g$ is the gravitational constant. Writing $\omega = 2\pi/T$ and $\kappa = 2\pi/L$, where $T$ is the wave period and $L$ the wavelength, the dispersion relation can also be expressed in terms of periods and wavelengths, namely

$$L = \frac{g}{2\pi}T^2.$$

Due to the dispersion relation, we thus have the following relation between the covariance function $R(\xi, \tau) = \mathsf{E}[W(x,t)W(x+\xi, t+\tau)]$ and the spectrum

$$(2.2) \qquad R(\xi, \tau) = \int_0^\infty \int_0^{2\pi} \cos\left(\omega\tau + \frac{\omega^2}{g}\xi\cos\theta\right)S(\omega, \theta)\,d\theta\,d\omega.$$

Due to this relation variances and covariances of the process and its derivatives can be expressed in terms of the spectral moments defined as

$$(2.3) \qquad \lambda_{ij} = \int_0^\infty \int_0^{2\pi} (\kappa(\omega)\cos\theta)^i\omega^j S(\omega, \theta)\,d\theta\,d\omega.$$

Note that $\mathrm{Var}(W(0,0)) = R(0,0) = \lambda_{00}$, $\mathrm{Var}(W_x(0,0)) = -R_{\xi\xi}(0,0) = \lambda_{20}$, $\mathrm{Var}(W_t(0,0)) = -R_{\tau\tau}(0,0) = \lambda_{02}$ and $\mathrm{Cov}(W_x(0,0), W_t(0,0)) = -R_{\xi\tau}(0,0) = \lambda_{11}$.

**3. Palm distributions and Rice's formula.** Let $N_Z^u([0,T])$ be the number of times that a process $Z(t)$ of one parameter takes the value $u$ in the interval $[0, T]$, and let $N_Z^u([0,T], E)$ be the number of times the process takes the value $u$ and at the same time a statement $E$ about the process, its derivatives or another process, is satisfied. For example $E$ could be the statement "$Z_t(t) \le z$," where $Z_t(t)$ is the derivative of $Z(t)$. A Palm measure is defined by the following ratio of intensities

$$(3.1) \qquad P_Z^u(E) = \frac{\mathsf{E}[N_Z^u([0,1], E)]}{\mathsf{E}[N_Z^u([0,1])]},$$

provided that the expected number of $u$-crossings, $\mathsf{E}[N_Z^u([0,1])]$, is finite. If the process $Z(t)$ is ergodic the intensities in (3.1) can be computed as sample averages and, with probability one,

$$P_Z^u(E) = \lim_{T\to\infty} \frac{N_Z^u([0,T], E)}{N_Z^u([0,T])}.$$

Thus the Palm distribution has the empirical interpretation that it is the long-term proportion of $u$-crossings by the process $Z$ for which $E$ is satisfied. Note that a sufficient condition for ergodicity for a Gaussian process is that it possesses a spectral density function; see [7], page 157.

The intensities of crossings in (3.1) can be computed by using Rice's formula, first studied by [8] and [15]. In its simplest form, when $Z$ is a zero-mean stationary Gaussian process having almost surely differentiable sample paths (realizations), it reads

$$(3.2) \qquad \mathsf{E}[N_Z^u([0,1])] = \frac{1}{\pi} \sqrt{\frac{\mathrm{Var}(Z_t(0))}{\mathrm{Var}(Z(0))}} \exp\left(-\frac{u^2}{2\,\mathrm{Var}(Z(0))}\right).$$

It can also be written on integral form, namely

$$
\begin{aligned}
(3.3) \qquad \mathsf{E}[N_Z^u([0,1])] &= \int_{-\infty}^{\infty} |z| f_{Z_t(0),Z(0)}(z,u)\,dz \\
&= \mathsf{E}[|Z_t(0)| \mid Z(0) = u] f_{Z(0)}(u),
\end{aligned}
$$

where $f_{Z_t(0),Z(0)}$ and $f_{Z(0)}$ are the densities of $(Z_t(0), Z(0))$ and $Z(0)$, respectively. There is a vast literature on generalized versions of Rice's formula; see, for example, [1, 3, 5, 9, 18]. It turns out that the formulation (3.3) extends in a natural way so that the intensity of crossings satisfying a statement $E$ can be computed. For example, if $Z(t)$ has almost surely differentiable sample paths and the event $E$ satisfies certain regularity conditions, then by [9], Lemma 7.5.2.

$$(3.4) \qquad \mathsf{E}[N_Z^u([0,1], E)] = \mathsf{E}[|Z_t(0)|\mathbf{1}_{\{E\}} \mid Z(0) = u] f_{Z(0)}(u),$$

where $\mathbf{1}_{\{\cdot\}}$ is an indicator function. Thus Rice's formula and generalizations of it can be used to compute the Palm distribution (3.1). However, some of the Palm distributions that we will consider include crossings by random vector fields and therefore slightly more general forms of Rice's formula are required. The exact version that will be used is due to [13], Theorem 9.6, which is a generalization of a theorem in [3].

REMARK 3.1. Notation like $\mathsf{E}[|Z_t(t)| \mid Z(t) = u] f_{Z(t)}(u)$ is used extensively throughout this paper. This should be interpreted in the sense of (3.3), that is, as an integral. In particular, this requires that the joint density of $Z_t(t)$ and $Z(t)$ exists, so that $(Z_t(t), Z(t))$ has a nondegenerate Gaussian distribution. For higher dimensions the notation should be understood in an analogous fashion.

REMARK 3.2. It should be pointed out that the Gaussian assumption is not necessary for Rice's formula (3.3) to hold, although the conditions for its validity in some cases can be quite intricate. However, if one is content to use (3.3) for almost every $u$, rather than for each specific $u$, these conditions simplify substantially. In [10] sufficient conditions for applying almost everywhere results to fixed levels are given.

**4. Wave velocities.** We will see that encountered distributions of wave characteristics are closely related to velocities of the individual apparent waves. Before defining velocities for these random waves we start with velocities of deterministic waves, and in particular, study the implications of the dispersion relation for the velocities of such waves.

Consider a deterministic cosine-wave $c(x,t)$ with angular wave frequency $\omega$ and wave number $\kappa$, namely

$$c(x,t) = A\cos(\omega t + \kappa x + \phi),$$

where $A$ is an amplitude and $\phi$ a phase. For fixed $t$ this is a wave in space and for fixed $x$, a wave in time. If both $t$ and $x$ are varying it is a wave traveling in the direction of the negative $x$-axis. The (phase) velocity of the wave can be expressed as

$$V = -\frac{L}{T} = -\frac{\omega}{\kappa} = -\frac{c_t(x,t)}{c_x(x,t)}.$$

Note that the convention here is that a wave travels in the direction of the positive $x$-axis if it has positive velocity and vice versa. Now for this deterministic wave, $\omega$ and $\kappa$ are linked by the dispersion relation (2.1) and consequently

$$(4.1) \qquad V = -\frac{\omega}{\kappa} = -\sqrt{\frac{g}{\kappa}} = -\sqrt{\frac{gL}{2\pi}}.$$

Thus for ocean waves the dispersion relation implies that longer waves are faster than shorter ones and, because the spatial slope of the wave at a zero-crossing is proportional to $\kappa$, it also implies that waves with a high spatial slope are slow.

Velocities for waves in irregular seas can be defined in a similar fashion, and has been studied by [2, 4, 12, 14]. In this case, however, one is interested in the velocities of the apparent waves rather than velocities of single harmonics. Recall that an apparent wave is the part of the sea record between two consecutive upcrossings of the still water level, and that the downcrossing in between is termed the center of the wave. Let $x_i > 0$ be the positions of centers of waves observed in $W(x,0)$, that is, along the $x$-axis at time zero. Due to the time variability of the sea surface, the centers change their positions with time and at $t = 0$ the velocity $V_i$ that the $i$th wave moves at can be evaluated by means of $V_i = -\frac{W_t(x_i,0)}{W_x(x_i,0)}$; see for example, [4]. The variability of the velocities $V_i$ will be described by the Palm distribution defined as follows

$$(4.2) \qquad F_V(v) = \frac{\mathsf{E}[\text{number of } x_i < 1 \text{ such that } V_i \leq v]}{\mathsf{E}[\text{number of } x_i < 1]}.$$

For the Gaussian sea this becomes

$$F_V(v) = \frac{1}{2}\left(1 + \frac{v - \bar{v}}{\sqrt{(v - \bar{v})^2 + \sigma^2/\lambda_{20}}}\right), \qquad \bar{v} = -\frac{\lambda_{11}}{\lambda_{20}},$$

where $\sigma^2 = \lambda_{02} - \frac{\lambda_{11}^2}{\lambda_{20}}$ and $\bar{v}$ is the average wave velocity; see [14].

REMARK 4.1. The velocity defined by $V(x,y,t) = -\frac{W_t(x,y,t)}{W_x(x,y,t)}$, $(x,y) \in \mathbb{R}^2$, can be seen as the local velocity field of a random, moving, surface. This velocity is defined at any point $(x,y) \in \mathbb{R}^2$ at any time $t$. For a thorough study of different velocity concepts; see [4].

**5. Explicit evaluation of Palm distributions.** Earlier it was shown how Palm distributions can be expressed by Rice's formula, given by (3.3) and (3.4). These formulas involve computation of multivariate normal expectations. Sometimes these have to be evaluated numerically, but in some cases explicit forms can be given. One of these explicit cases will be investigated next in some detail.

In the following lemma, a stationary Gaussian vector valued process $(Z(t), Z_t(t), Y(t))$ will be considered. In order to simplify notation, introduce the random variables $Y^u$ and $Z_t^u$ having joint distribution function $F_{Y^u, Z_t^u}(y, z)$ defined by the measure (3.1), namely

$$(5.1) \qquad F_{Y^u, Z_t^u}(y, z) = P(Y^u \le y, Z_t^u \le z) = P_Z^u(Y(t) \le y, Z_t(t) \le z).$$

These random variables are a special case of the Slepian model process; see [11] for a detailed presentation, and one may think of them as $Y(t)$ and $Z_t(t)$ observed at a randomly chosen $u$-crossing of the process $Z(t)$. The next lemma gives an explicit representation of these variables.

LEMMA 5.1. *Suppose that $(Z(t), Z_t(t), Y(t))$ is a stationary Gaussian vector valued process such that the density of $(Z(t), Z_t(t), Y(t))$ exists. Assume further that the sample paths of $Y(t)$ and $Z(t)$ are a.s. continuous and a.s. continuously differentiable, respectively. Let $m_Z = \mathsf{E}(Z(t))$, $m_Y = \mathsf{E}(Y(t))$ and let the covariance matrix $\Sigma$ of the vector $(Z(t), Z_t(t), Y(t))$ be given by*

$$(5.2) \qquad \Sigma = \begin{pmatrix} \sigma_Z^2 & 0 & \sigma_Z \sigma_Y \rho_{ZY} \\ 0 & \sigma_{Z_t}^2 & \sigma_{Z_t} \sigma_Y \rho_{Z_t Y} \\ \sigma_Z \sigma_Y \rho_{ZY} & \sigma_{Z_t} \sigma_Y \rho_{Z_t Y} & \sigma_Y^2 \end{pmatrix}.$$

*Then the variables $Y^u, Z_t^u$, with distribution defined by (5.1), have the following distributional representations:*

$$(5.3) \qquad Z_t^u = \sigma_{Z_t} R, \qquad Y^u = m_u + \sigma_Y(\rho_{Z_t Y} R + \sqrt{1 - \rho_{ZY}^2 - \rho_{Z_t Y}^2} U),$$

*where $U$ is standard normal and $R$ is independent of $U$, with a double Rayleigh distribution, that is, probability density $f(r) = \frac{|r|}{2} e^{-r^2/2}$. Here $m_u = m_Y + \frac{\sigma_Y}{\sigma_Z} \rho_{ZY}(u - m_Z)$ and $\rho_{ZY}$, $\rho_{Z_tY}$ are correlations between $Y(0)$ and $Z(0)$, $Z_t(0)$, respectively.*

PROOF.  The proof is easy using the definition of Palm probability, Rice's formula (3.2) and (3.4) and the Gaussian regression formulas.  □

This lemma thus shows that the derivative $Z_t$ at crossings has a Rayleigh distribution as opposed to the derivative at any point, where it has a Gaussian distribution. Furthermore, the distribution of another Gaussian process $Y$, correlated with $Z(t)$, can be represented as the sum of a Rayleigh and Gaussian variable at the crossing point. The corresponding distributions, but only observed at up- or downcrossings of $Z(t)$, are easily obtained by using the representation from Lemma 5.1. These distributions are summarized in the following lemma.

LEMMA 5.2.  *Let $Z_t^u \mid Z_t^u > 0$ denote the random variable obtained by conditioning $Z_t^u$ on $Z_t^u > 0$ and let $Y^u \mid Z_t^u > 0$, $Z_t^u \mid Z_t^u < 0$ and $Y^u \mid Z_t^u < 0$ be defined analogously. Using the same notation as in Lemma 5.1, we then have the following distributional representations:*

$$Z_t^u \mid Z_t^u > 0 = \sigma_{Z_t} R^+,$$

$$Y^u \mid Z_t^u > 0 = m_u + \sigma_Y(\rho_{Z_tY} R^+ + \sqrt{1 - \rho_{ZY}^2 - \rho_{Z_tY}^2}\, U)$$

*and*

$$Z_t^u \mid Z_t^u < 0 = \sigma_{Z_t} R^-,$$

$$Y^u \mid Z_t^u < 0 = m_u + \sigma_Y(\rho_{Z_tY} R^- + \sqrt{1 - \rho_{ZY}^2 - \rho_{Z_tY}^2}\, U),$$

*where $R^+$ and $R^-$, both independent of $U$, have densities $f_{R^+}(r) = re^{-r^2/2}$, $r > 0$, and $f_{R^-}(r) = -re^{-r^2/2}$, $r < 0$, respectively.*

This result will be used to obtain explicit representations of the distributions of the slope observed at wave centers and encountered wave centers. When evaluating such distributions, the following lemma that can be proved by straightforward calculations is useful.

LEMMA 5.3.  *Let $U$ and $R^+$ be independent standard Gaussian and Rayleigh distributed variables. Then, for $a > 0$ and any $x$,*

$$P(aU + bR^+ > x) = \Phi\left(-\frac{x}{a}\right) + e^{-x^2/(2\sigma^2)} \cdot \frac{b}{\sigma} \Phi\left(\frac{x}{\sigma}\frac{b}{a}\right),$$

*where $\sigma = \sqrt{a^2 + b^2}$ and $\Phi$ is the distribution function of a standard Gaussian variable.*

**6. Distribution of slope.** Consider the sea surface $W(x,t)$. In this section the distribution of the slope $W_x(x,t)$ will be computed, but only for values of $x$ and $t$ chosen in a careful manner. The cases considered are:

1. *Spatial case.* Choose $x$ to be points where the process $W(x,0)$ has a downcrossing of the zero level. This means that only values of $x$ satisfying $W(x,0) = 0$ and $W_x(x,0) < 0$ are chosen. In other words, $x$ are locations of the wave centers in $W(\cdot,0)$. These values of $x$ and $t = 0$ will lead to a Palm distribution of wave slopes observed at centers of waves.

2. *Encountered case.* Consider a ship sailing with constant velocity $v$ on the ocean. If it at time zero is at position $x = 0$, then the sea elevation at the center of gravity at time $t$ is $W(vt,t)$. Choose $t$ so that $W(vt,t) = 0$, $W_x(vt,t) < 0$ and $\frac{\partial}{\partial t}W(vt,t) < 0$, then $t$ can be interpreted as a time when the ship is overtaken by a wave center. Values of $t$ chosen in this manner, and $x = vt$, lead to a Palm distribution of wave slopes observed at overtaking centers of waves.

REMARK 6.1. By slope is here meant the water surface wave slope. This should not be confused with the term velocity slope sometimes used in the marine sciences.

6.1. *Spatial distribution of slope.* To derive the distribution of the slope observed at centers of waves, we start by defining it by means of a Palm distribution. Therefore, let $x_i \geq 0$ be the locations of the centers of waves, that is, locations of the zero-downcrossings, in $W(x,0)$. Each of these wave centers can be associated with a slope $W_x(x_i,0)$ and the following Palm distribution for the slope at centers of waves, denoted by $F_{W_x}^{\text{space}}(w)$, can be defined as

$$\begin{aligned}
(6.1) \quad & F_{W_x}^{\text{space}}(w) \\
&= \frac{\mathsf{E}[\text{number of } x_i \leq 1 \text{ such that } W_x(x_i,0) \leq w]}{\mathsf{E}[\text{number of } x_i \leq 1]} \\
&= \frac{\mathsf{E}[\#\{x \in [0,1]; W(x,0) = 0, W_x(x,0) < 0, W_x(x,0) \leq w\}]}{\mathsf{E}[\#\{x \in [0,1]; W(x,0) = 0, W_x(x,0) < 0\}]},
\end{aligned}$$

where the second equality holds because $x_i$ is a point of downcrossing of $W(x,0)$. Note that by ergodicity this distribution can be interpreted as

$$F_{W_x}^{\text{space}}(w) = \lim_{X \to \infty} \frac{\text{number of } x_i < X \text{ such that } W_x(x_i,0) \leq w}{\text{number of } x_i < X},$$

that is, as the proportion of wave centers with associated slope less or equal to $w$ in an infinitely long realization of the process.

The following theorem gives an expression for the Palm distribution in this case, and also states how it can be evaluated.

THEOREM 6.1.   *If $W(x,0)$ is a stationary Gaussian process having a.s. continuously differentiable sample paths then*

$$(6.2) \quad F_{W_x}^{\text{space}}(w) = \frac{\mathsf{E}[W_x(0,0)^- \mathbf{1}_{\{W_x(0,0) \leq w\}} \mid W(0,0) = 0] f_{W(0,0)}(0)}{\mathsf{E}[W_x(0,0)^- \mid W(0,0) = 0] f_{W(0,0)}(0)},$$

*where $x^- = \max(-x, 0)$. Moreover*

$$(6.3) \qquad\qquad F_{W_x}^{\text{space}}(w) = P(\sqrt{\lambda_{20}} R^- \leq w),$$

*where $R^-$ is a random variable having density $f_{R^-}(r) = -re^{-r^2/2}, \ r < 0$.*

PROOF.   The first statement of the proof follows by applying the generalized Rice's formula (3.4) to the numerator and denominator of (6.1). For the numerator, (3.4) is used with $E = $ "$W_x(x,0) < 0, W_x(x,0) \leq w$" and for the denominator $E = $ "$W_x(x,0) < 0$."

To prove (6.3), Lemma 5.1 is used with $Z = W(x,0)$. Then $Z_x = W_x(x,0)$, and the Slepian variable $Z_x^0$ is the derivative $W_x(x,0)$ observed at zero-crossings of the process $Z = W(x,0)$. Using the variable $Z_x^0$, the distribution for the slope $W_x(x,0)$ observed at zero-downcrossings can be expressed as

$$F_{W_x}^{\text{space}}(w) = P(Z_x^0 \leq w \mid Z_x^0 < 0),$$

that is, as the distribution of the slope observed at zero-crossings conditional on a negative slope at the crossing. Now (6.3) follows by Lemma 5.2.   □

6.2. *Encountered distribution of slope.*   Consider a vessel sailing along the $x$-axis in the positive direction with constant speed $v$, having its center of gravity at time zero at $x = 0$. Disregarding the interaction between the ship and the waves, the sea elevation at the center of gravity of the ship is given by $Z(t) = W(vt, t)$. In the following, we will refer to $Z(t)$ as the *encountered sea.* Now let $t_i \geq 0$ be the times when the vessel is overtaken by a wave center. Characteristic for such times is that the center of gravity is passing through the still water level, $Z(t_i) = 0$, the slope of the wave is negative, $W_x(vt_i, t_i) < 0$ and the encountered sea $Z(t)$ has an upcrossing $Z_t(t_i) > 0$, meaning that the center of wave is overtaking the ship. Because each $t_i$ can be associated with the slope $W_x(vt_i, t_i)$ of the $i$th overtaking wave, the following Palm distribution of the slope at overtaking waves, $F_{W_x}^{\text{enc}}(w)$, say, can be formed

$$
\begin{aligned}
& F_{W_x}^{\text{enc}}(w) \\
(6.4) \quad & = \frac{\mathsf{E}[\text{number of } t_i \leq 1 \text{ such that } W_x(vt_i, t_i) \leq w]}{\mathsf{E}[\text{number of } t_i \leq 1]} \\
& = \frac{\mathsf{E}[\#\{t \in [0,1]; Z(t) = 0, Z_t(t) > 0, W_x(vt, t) < 0, W_x(vt, t) \leq w\}]}{\mathsf{E}[\#\{t \in [0,1]; Z(t) = 0, Z_t(t) > 0, W_x(vt, t) < 0\}]},
\end{aligned}
$$

where the second equality is due to the definition of $t_i$. As before, the Palm distribution can be interpreted in the following frequency fashion:

$$F_{W_x}^{\text{enc}}(w) = \lim_{T \to \infty} \frac{\text{number of } t_i < T \text{ such that } W_x(vt_i, t_i) \leq w}{\text{number of } t_i < T}.$$

REMARK 6.2. Despite the similarities of expressions (6.1) and (6.4), there is an important difference arising from the fact that crossings from two different processes are counted, namely crossings from $W(x, 0)$ in the spatial case and from $W(vt, t)$ in the encountered case. From a statistical point of view this means sampling from two different populations.

The following theorem, concerning the Palm distribution of slope for encountered waves, is an analogue to Theorem 6.1.

THEOREM 6.2.    Let $W(x, t)$ be a stationary Gaussian process having a.s. continuously differentiable sample paths. Then $F_{W_x}^{\text{enc}}(w)$ can be expressed as

$$(6.5) \quad F_{W_x}^{\text{enc}}(w) = \frac{\mathsf{E}[Z_t(0)^+ \mathbf{1}_{\{W_x(0,0)<0, W_x(0,0)\leq w\}} \mid Z(0)=0] f_{Z(0)}(0)}{\mathsf{E}[Z_t(0)^+ \mathbf{1}_{\{W_x(0,0)<0\}} \mid Z(0)=0] f_{Z(0)}(0)},$$

where $x^+ = \max(x, 0)$. Moreover, $F_{W_x}^{\text{enc}}(w) = 1$ if $w \geq 0$ and for $w < 0$

$$(6.6) \quad \begin{aligned} F_{W_x}^{\text{enc}}(w) = \frac{2}{1 - \rho_{Z_t Y}} \bigg( &\Phi\bigg(\frac{w}{\sqrt{\lambda_{20}(1 - \rho_{Z_t Y}^2)}}\bigg) \\ &- \rho_{Z_t Y} e^{-w^2/(2\lambda_{20})} \Phi\bigg(\frac{\rho_{Z_t Y} w}{\sqrt{\lambda_{20}(1 - \rho_{Z_t Y}^2)}}\bigg)\bigg), \end{aligned}$$

where $\Phi$ is the distribution function of a standard normal variable and

$$(6.7) \quad \rho_{Z_t Y} = \frac{v\lambda_{20} + \lambda_{11}}{\sqrt{\lambda_{20}(v^2\lambda_{20} + 2v\lambda_{11} + \lambda_{02})}}.$$

PROOF.    Expression (6.5) follows, in a similar fashion as in the proof of Theorem 6.1, by applying the generalized Rice's formula (3.4) to the numerator and denominator in (6.4).

To prove the second statement, identify the process $Z(t)$ of Lemma 5.1 with the encountered process $W(vt, t)$ and let $Y(t) = W_x(vt, t)$. With this choice of variables the covariance matrix, expressed in terms of spectral moments defined by (2.3), becomes

$$(6.8) \quad \Sigma = \begin{pmatrix} \lambda_{00} & 0 & 0 \\ 0 & v^2\lambda_{20} + 2v\lambda_{11} + \lambda_{02} & v\lambda_{20} + \lambda_{11} \\ 0 & v\lambda_{20} + \lambda_{11} & \lambda_{20} \end{pmatrix}.$$

In this setting $Z_t^0$ is the time derivative of the encountered process observed at its zero-crossings, and $Y^0$ is the spatial slope observed at zero-crossings of the encountered process. The distribution in (6.5), however, is the distribution of the spatial slope, but only at upcrossings of $Z(t)$ such that the spatial slope is negative. Using the variables $Z_t^0$ and $Y^00$ this can be expressed as

$$F_{W_x}^{\text{enc}}(w) = P(Y^0 \leq w \mid Z_t^0 > 0, Y^0 < 0)$$

$$= \begin{cases} \dfrac{P(Y^0 \leq w \mid Z_t^0 > 0)}{P(Y^0 < 0 \mid Z_t^0 > 0)}, & w < 0, \\ 1, & w \geq 0. \end{cases}$$

Use of Lemma 5.2 now gives the following representation

$$F_{W_x}^{\text{enc}}(w) = \begin{cases} \dfrac{P(\sqrt{\lambda_{20}}(\rho_{Z_tY}R^+ + \sqrt{1-\rho_{Z_tY}^2}U) \leq w)}{P(\sqrt{\lambda_{20}}(\rho_{Z_tY}R^+ + \sqrt{1-\rho_{Z_tY}^2}U) < 0)}, & w < 0, \\ 1, & w \geq 0, \end{cases}$$

with $\rho_{Z_tY}$ given by (6.7). Thus the slope observed at an encountered overtaking center of wave can be evaluated by computing probabilities for a sum of a Rayleigh and a Gaussian random variable. Thus, (6.6) holds by Lemma 5.3. $\square$

The interpretation of (6.5) is that it is the distribution of the slope of the waves that overtake a ship sailing on the ocean. If a wave overtake the ship, then it must necessarily be traveling at a higher speed than the ship. Thus intuition suggests that in some way the distribution (6.5) should be related to the relative velocities of the individual waves and the ship. The exact relation is given by the following corollary, enabling us to interpret the encountered distribution in terms of physical quantities.

COROLLARY 6.1. *Let* $V = -W_t(0,0)/W_x(0,0)$ *be the local velocity at* $(x,t) = (0,0)$ *and let* $v$ *denote the velocity of the ship. Then*

$$(6.9) \quad \begin{aligned} &F_{W_x}^{\text{enc}}(w) \\ &= \frac{\mathsf{E}[(V-v)^+ W_x(0,0)^- \mathbf{1}_{\{W_x(0,0) \leq w\}} \mid W(0,0) = 0] f_{W(0,0)}(0)}{\mathsf{E}[(V-v)^+ W_x(0,0)^- \mid W(0,0) = 0] f_{W(0,0)}(0)}. \end{aligned}$$

PROOF. The assertion of the corollary follows easily from (6.5) due to the following equality

$$\begin{aligned} (6.10) \quad Z_t(0)^+ \mathbf{1}_{\{W_x(0,0)<0\}} &= (vW_x(0,0) + W_t(0,0))^+ \mathbf{1}_{\{W_x(0,0)<0\}} \\ &= (W_x(0,0)(v-V))^+ \mathbf{1}_{\{W_x(0,0)<0\}} \\ &= W_x(0,0)^-(v-V)^- \\ &= W_x(0,0)^-(V-v)^+, \end{aligned}$$

and the fact that $Z(0) = W(0,0)$.  □

The encountered distribution given by (6.5) or (6.9) differs from the corresponding spatial distribution (6.2) only by the term $(V - v)^+$ that enters into the expectation. Thus the spatial measure for slope can be transformed to the encountered measure for slope by entering the velocity term $(V - v)^+$. The physical interpretation of this is that in the encountered approach not every center of wave is sampled but only those with velocity greater than the ship. The velocity factor is thus a consequence of the fact that the sea surface is observed from a moving observer, that is, the ship. In physics, such a phenomenon would be termed a *Doppler shift*. In the following, we will see that the Doppler shift transformation of the distribution is not limited to the distribution of the spatial slope. In fact, by the way the measures are defined, it extends to all wave characteristics of overtaking encountered waves.

## 7. Distribution of waveheight and half-wavelength.
Let $x_1$ be the location of a wave center, that is, a downcrossing of the level zero, in $W(x,t)$ for some time $t$. Then $x_1$ is characterized by the fact that $W(x_1, t) = 0$ and $W_x(x_1, t) < 0$. Further, let $x_2$ be the distance from the downcrossing $x_1$ to the closest local maximum before $x_1$ so that $W_x(x_1 - x_2, t) = 0$, $W_{xx}(x_1 - x_2, t) < 0$ and $W_x(x_1 - s, t) \leq 0, \forall s \in (0, x_2)$. The latter condition guarantees that this really is the closest maximum before $x_1$. Similarly define $x_3$ to be the distance from $x_1$ to the closest local minimum after $x_1$. Hence $W_x(x_1 + x_3, t) = 0$, $W_{xx}(x_1 + x_3, t) > 0$ and $W_x(x_1 - s, t) \leq 0, \forall s \in (0, x_3)$. The definitions of $x_2$ and $x_3$ are shown in Figure 1. In the same figure, it is also indicated that the notation $H_2 = W(x_1 - x_2, t)$ and $H_3 = W(x_1 + x_3, t)$ will be used.

The task is to obtain an expression for the joint Palm distribution of the distances from the wave center to the two closest local extrema and their heights $(x_2, x_3, H_2, H_3)$. As before two cases, corresponding to choosing wave centers in two different ways, are considered. The first is a spatial case where the Palm distribution is derived for wave centers in $W(x, 0)$ and the second is an encountered case where the corresponding distribution is derived for wave centers that overtake a ship.

By dividing the waveheight by the half-wavelength a measure of the steepness at the wave center is obtained, hence motivating the study of these quantities from a ship stability perspective.

### 7.1. *Spatial distribution of heights and distances.*
In this section, an expression for the joint Palm distribution of $(x_2, x_3, H_2, H_3)$ at downcrossings in space is derived. To do so, we first define the appropriate Palm distribution and then express it in terms of level crossings of the process $W(x) = W(x, 0)$.

Let $(x_1)_i \geq 0$ be positions of downcrossings in $W(x)$. With each downcrossing associate distances $(x_2, x_3)_i$, say, to the closest local extrema and form the following Palm distribution:

$$F_{x_2, x_3, H_2, H_3}^{\text{space}}(r, s, u, w)$$
$$= \frac{\mathsf{E}[\text{number of } (x_1)_i \leq 1 \text{ such that } (x_2, x_3, H_2, H_3)_i \leq (r, s, u, w)]}{\mathsf{E}[\text{number of } (x_1)_i \leq 1]}.$$

In order to express the Palm distribution in terms of level crossings of the process $W$, note that for each crossing $(x_1)_i$ satisfying $(x_2, x_3, H_2, H_3)_i \leq (r, s, u, w)$ the triple $(x_1, x_2, x_3)$ is characterized by the following properties:

(A1)      $(x_1, x_2, x_3) \in B$,  where the set $B$ is defined by

$$B = \{x \in \mathbb{R}^3; 0 \leq x_1 \leq 1, 0 < x_2 \leq r, 0 < x_3 \leq s\},$$

(A2)      $W(x_1) = W_x(x_1 - x_2) = W_x(x_1 + x_3) = 0,$

(A3)      $\begin{aligned} &W_{xx}(x_1 - x_2) < 0, \\ &W_{xx}(x_1 + x_3) > 0, \end{aligned}$

$W_x(s) \leq 0 \qquad \forall s \in \Gamma_x,$  where $\Gamma_x = \{s \in \mathbb{R}; x_1 - x_2 \leq s \leq x_1 + x_3\},$

(A4)      $W(x_1 - x_2) \leq u, \qquad W(x_1 + x_3) \leq w.$

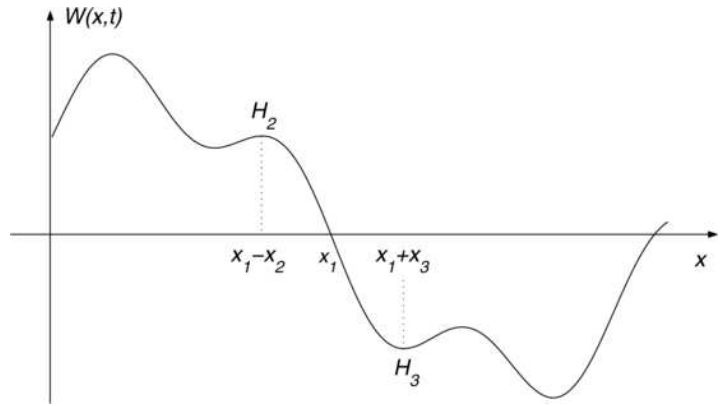

FIG. 1.   *Definition of waveheights and distances. If $x_1$ is the location of a downcrossing in $W(x, t)$ for fixed $t$, then $x_2$ is defined to be the distance to the first local maximum before $x_1$ and similarly $x_3$ is the distance to the first local minimum after $x_1$. $H_2$ is defined as the sea elevation $W(x_1 - x_2, t)$ and $H_3$ as the sea elevation $W(x_1 + x_3, t)$.*

Defining

$$(7.1) \qquad \xi(x) = \xi(x_1, x_2, x_3) = (W(x_1), W_x(x_1 - x_2), W_x(x_1 + x_3))$$

it thus holds that

$$(7.2) \qquad \begin{aligned} & F^{\text{space}}_{x_2, x_3, H_2, H_3}(r, s, u, w) \\ & \quad = \frac{\mathsf{E}[\#\{x \in B; \xi(x) = 0, (\text{A3}) \text{ and } (\text{A4}) \text{ satisfied}\}]}{\mathsf{E}[\#\{x_1 \in [0, 1]; W(x_1) = 0, W_x(x_1) < 0\}]}. \end{aligned}$$

The formulation $(7.2)$ is suitable when it comes to computation of the expectations by means of Rice's formula. In particular, the following theorem states how the expectation in the numerator can be computed.

THEOREM 7.1. *Let $W : \mathbb{R} \to \mathbb{R}$ be a stationary Gaussian process having a.s. twice continuously differentiable sample paths. Assume that the spectrum of $W$ has a continuous component and that the variance of the number of zeros $W_x(y) = 0$, $y \in [0, 1]$ is finite. From $W$ define a process $\xi$ by the relation $(7.1)$ and a set $B$ as in* (A1). *Furthermore, define a vector valued process $Y = (W, W_x, W_{xx})$, a set $\Gamma_x$ as in* (A3) *and let $g(Y, x)$ be an indicator function defined by*

$$\begin{aligned} g(Y, x) = {} & \mathbf{1}_{\{W_x(s) \le 0, \forall s \in \Gamma_x\}} \mathbf{1}_{\{W(x_1 - x_2) \le u\}} \mathbf{1}_{\{W(x_1 + x_3) \le w\}} \\ & \times \mathbf{1}_{\{W_{xx}(x_1 - x_2) < 0\}} \mathbf{1}_{\{W_{xx}(x_1 + x_3) > 0\}}. \end{aligned}$$

*Then, writing $N_0^{\xi}(B, g) = \#\{x \in B; \xi(x) = 0, g(Y, x) = 1\}$ and $\det(\xi_x(x))$ for the Jacobian determinant of $\xi(x)$,*

$$\mathsf{E}[N_0^{\xi}(B, g)] = \int_B \mathsf{E}[|\det(\xi_x(x))| g(Y, x) \mid \xi(x) = 0] f_{\xi(x)}(0) \, dx,$$

*where both members are finite.*

PROOF. Let $B_{\delta} = \{x \in \mathbb{R}^3; 0 \le x_1 \le 1, 0 < x_2 < \frac{\delta}{2}, 0 < x_3 < \frac{\delta}{2}\}$ and $B_{-\delta} = B \backslash B_{\delta}$. Then

$$(7.3) \qquad \sum_{x \in B \cap \xi^{-1}(0)} g(Y, x) = \sum_{x \in B_{\delta} \cap \xi^{-1}(0)} g(Y, x) + \sum_{x \in B_{-\delta} \cap \xi^{-1}(0)} g(Y, x).$$

On the set $B_{-\delta}$ the distribution of $\xi$ is nondegenerate, so that Theorem 9.6. of [13] can be used. However, that theorem is only valid for a continuous function $g$. Therefore, let $\{k_n\}$ and $\{h_n\}$ be sequences of continuous, monotone functions such that $k_n(x) = 1$ if $x \le 0$, $k_n(x) = 0$ if $x > 1/n$, $h_n(x) = 1$ if $x \le -1/n$ and $h_n(x) = 0$ if $x > 0$. It is easy to verify that $k_n(x) \to \mathbf{1}_{(-\infty, 0]}(x)$ and $h_n(x) \to \mathbf{1}_{(-\infty,0)}(x)$ as $n \to \infty$. Because $k_n(x)$ and $h_n(x)$ are continuous for each $n$ the theorem of [13] can be applied with

$$g_n(Y,x) = k_n\left(\sup_{s \in \Gamma_x} W_x(s)\right) \cdot k_n(W(x_1 - x_2) - u) \cdot k_n(W(x_1 + x_3) - w)$$
$$\times h_n(W_{xx}(x_1 - x_2)) \cdot h_n(-W_{xx}(x_1 + x_3)),$$

so that

$$(7.4) \quad \begin{aligned} &\mathsf{E}\left[\sum_{x \in B_{-\delta} \cap \xi^{-1}(0)} g_n(Y,x)\right] \\ &= \int_{B_{-\delta}} \mathsf{E}[|\det(\xi_x(x))|g_n(Y,x) \mid \xi(x) = 0]f_{\xi(x)}(0) \, dx, \end{aligned}$$

where both members are finite. By dominated convergence arguments as $n \to \infty$, (7.4) holds with $g_n$ replaced by $g$.

On $B_\delta$, however, the distribution of $\xi$ will degenerate, so the same theorem cannot be used in this case. Instead, we show that the assumption of finite variance of the number of zeros $W_x(y) = 0$, $y \in [0,1]$, will force the expected value of the first term on the right-hand side in (7.3) to go to zero as $\delta$ approaches zero. To see that this is the case, first note that $g$ is bounded by 1 so that

$$0 \le \mathsf{E}\left[\sum_{x \in B \cap \xi^{-1}(0)} g(Y,x)\right] - \mathsf{E}\left[\sum_{x \in B_{-\delta} \cap \xi^{-1}(0)} g(Y,x)\right] \le \mathsf{E}[N_0^\xi(B_\delta)],$$

where $N_0^\xi(B_\delta) = \#\{x \in B_\delta; \xi(x) = 0\}$. Because $W$ is assumed to be an a.s. differentiable Gaussian process, it follows that a.s. there are only finitely many solutions to $W(y) = 0$ in any finite interval. Thus, with probability one, $\{y \in [0,1]; W(y) = 0\} = \{y_1, \ldots, y_K\}$, where $K = N_0^W([0,1]) = \#\{y \in [0,1]; W(y) = 0\}$. Using this fact and the stationarity of $W$, the following inequality holds for $0 < \delta < 1$:

$$\mathsf{E}[N_0^\xi(B_\delta)] = \mathsf{E}\left[\sum_{i=1}^K \#\{(x_2,x_3) \in [0,\delta/2]^2; W_x(y_i - x_2) = W_x(y_i + x_3) = 0\}\right]$$

$$= \sum_{k=1}^\infty \sum_{i=1}^k \mathsf{E}[\#\{(x_2,x_3) \in (0,\delta/2]^2;$$
$$W_x(y_i - x_2) = W_x(y_i + x_3) = 0\}]P(K = k)$$

$$= \sum_{k=1}^\infty k\mathsf{E}[\#\{(x_2,x_3) \in (0,\delta/2]^2; W_x(-x_2) = W_x(x_3) = 0\}]P(K = k)$$

$$= \mathsf{E}[\#\{(x_2,x_3) \in (0,\delta/2]^2; W_x(-x_2) = W_x(x_3) = 0\}] \cdot \mathsf{E}[N_0^W([0,1])]$$

$$\leq \mathsf{E}[\#\{(s,t) \in [0,1];$$
$$0 < |t-s| < \delta, W_x(s,0) = W_x(t,0) = 0\}] \cdot \mathsf{E}[N_0^W([0,1])].$$

Writing $N_\delta = \#\{(s,t) \in [0,1]; 0 < |t-s| < \delta, W_x(s,0) = W_x(t,0) = 0\}$ and observing that

$$N_0^{W_x}([0,1])(N_0^{W_x}([0,1]) - 1)$$
$$= \#\{(s,t) \in [0,1]; s \neq t, W_x(s,0) = W_x(t,0) = 0\},$$

where $N_0^{W_x}([0,1]) = \#\{y \in [0,1]; W_x(y) = 0\}$, one can by using the formula for the second factorial moment, given in [7], deduce that

$$\mathsf{E}[N_\delta] = \int_{\{(t,s) \in [0,1]^2, 0 < |t-s| < \delta\}} h(s,t)\, ds\, dt,$$

where

$$h(s,t) = \mathsf{E}[|W_{xx}(s)W_{xx}(t)| \mid W_x(s) = W_x(t) = 0] f_{W_x(s),W_x(t)}(0,0).$$

Because the assumption of finite variance of the number of zeros $N_0^{W_x}([0,1])$ implies that $h(s,t)$ is integrable over the rectangle $[0,1]^2$ it follows, by absolute continuity, that $\lim_{\delta \to 0} \mathsf{E}[N_\delta] = 0$. Thus

$$\mathsf{E}\left[\sum_{x \in B \cap \xi^{-1}(0)} g(Y,x)\right]$$
$$= \lim_{\delta \to 0} \int_{B_{-\delta}} \mathsf{E}[|\det(\xi_x(x))|g(Y,x) \mid \xi(x) = 0] f_{\xi(x)}(0)\, dx,$$

so that the final result follows by monotone convergence. □

REMARK 7.1. In [13] the domain of definition of the function $g$ is $\mathcal{C}(B, \mathbb{R}) \times B$, where $\mathcal{C}(B, \mathbb{R})$ is the class of continuous functions on $B$. However, the arguments of that proof holds also in the case when the domain of definition of $g$ is $\mathcal{C}(B, \mathbb{R}^n) \times B$, where $\mathcal{C}(B, \mathbb{R}^n)$ is the class of continuous functions on $B$ taking values in $\mathbb{R}^n$. The spaces $\mathcal{C}(B, \mathbb{R})$ and $\mathcal{C}(B, \mathbb{R}^n)$ are, in both cases, equipped with a uniform convergence norm.

REMARK 7.2. Sufficient conditions for finite variance of the number of zeros $W_x(y) = 0$, $y \in [0,1]$, can be found in [7], page 209.

Using Theorem 7.1, we may state the main theorem for the Palm distribution of $(x_2, x_3, H_2, H_3)$ in the spatial case.

THEOREM 7.2. *Let the assumptions and notation of Theorem 7.1 hold and define*

$$(7.5) \quad \alpha(x) = W_x(x_1)^- W_{xx}(x_1 - x_2)^- W_{xx}(x_1 + x_3)^+ \mathbf{1}_{\{W_x(s) \leq 0, \forall s \in \Gamma_x\}}.$$

*Moreover, define a function $h_s(x,y) = h_s(x_1, x_2, x_3, y_1, y_2)$ by*

$$h_s(x,y) = \mathsf{E}[\alpha(x) \mid \xi(x) = 0, W(x_1 - x_2) = y_1, W(x_1 + x_3) = y_2]$$

$$\times f_{\xi(x), W(x_1 - x_2), W(x_1 + x_3)}(0, y_1, y_2).$$

*Then the Palm distribution $F^{\mathrm{space}}_{x_2, x_3, H_2, H_3}$ is given by*

$$(7.6) \quad \begin{aligned} &F^{\mathrm{space}}_{x_2, x_3, H_2, H_3}(r, s, u, w) \\ &= \frac{\int_0^r \int_0^s \int_0^u \int_{-\infty}^w h_s(0, x_2, x_3, y_1, y_2)\, dy_2\, dy_1\, dx_3\, dx_2}{\mathsf{E}[W_x(0)^- \mid W(0) = 0] f_{W(0)}(0)}, \end{aligned}$$

*having the density*

$$(7.7) \quad f^{\mathrm{space}}_{x_2, x_3, H_2, H_3}(r, s, u, w) = \frac{h_s(0, r, s, u, w)}{\mathsf{E}[W_x(0)^- \mid W(0) = 0] f_{W(0)}(0)},$$

*where $r, s, u > 0$ and $w < 0$.*

PROOF. From the definition of the Palm distribution (7.2) and Theorem 7.1, the numerator in the statement (7.6) of the theorem follows due to stationarity in the first coordinate $x_1$. The denominator follows by straightforward application of Rice's formula; see the proof of Theorem 6.1. Moreover, it follows from (7.6) that $F^{\mathrm{space}}$ is absolutely continuous so that the density exists and is given by (7.7). □

7.2. *Encountered distribution of heights and distances.* The corresponding Palm distribution for overtaking encountered waves will now be considered, that is the distribution of distances and heights at times $t_i \geq 0$ when a ship sailing with constant velocity $v$ is overtaken by a wave center. To obtain a Palm distribution each crossing $t_i$ is associated with some distances denoted by $(x_2, x_3)_i$. In this case $(x_2, x_3)_i$ are the distances from $vt_i$ to the closest local maximum and minimum before and after $vt_i$, respectively, in the wave profile $W(x, t_i)$. Following now familiar steps, the following Palm distribution can be formed

$$\begin{aligned} &F^{\mathrm{enc}}_{x_2, x_3, H_2, H_3}(r, s, u, w) \\ &= \frac{\mathsf{E}[\text{number of } t_i \leq 1 \text{ such that } (x_2, x_3, H_2, H_3)_i \leq (r, s, u, w)]}{\mathsf{E}[\text{number of } t_i \leq 1]}. \end{aligned}$$

Similarly to the previous case, the Palm distribution can be expressed in terms of level crossings of the process $W$. More precisely, each time $t_i$ satisfying $(x_2, x_3, H_2, H_3)_i \le (r, s, u, w)$ has the following characteristics

(C1)
$$(t, x_2, x_3) \in B, \text{ where the set } B \text{ is defined by}$$
$$B = \{(t, x_2, x_3) \in \mathbb{R}^3; 0 \le t \le 1, 0 < x_2 \le r, 0 < x_3 \le s\},$$

(C2)
$$W(vt, t) = W_x(vt - x_2, t) = W_x(vt + x_3, t) = 0,$$

$$W_{xx}(vt - x_2, t) < 0,$$

(C3)
$$W_{xx}(vt + x_3, t) > 0,$$

$$W_x(s, t) \le 0, \qquad \forall s \in \Gamma_{tx}, \text{ where } \Gamma_{tx} = \{s \in \mathbb{R}; vt - x_2 \le s \le vt + x_3\}$$

(C4)
$$W(vt - x_2, t) \le u, \qquad W(vt + x_3, t) \le w,$$

(C5)
$$\frac{\partial}{\partial t} W(vt, t) = v W_x(vt, t) + W_t(vt, t) > 0.$$

Using these facts and by defining

(7.8)    $\eta(x) = \eta(t, x_2, x_3) = (W(vt, t), W_x(vt - x_2, t), W_x(vt + x_3, t))$

the Palm distribution can be written in the following form, suitable for our purposes:

(7.9)
$$F^{\text{enc}}_{x_2, x_3, H_2, H_3}(r, s, u, w)$$
$$= \frac{\mathsf{E}[\#\{x \in B; \eta(x) = 0, (\text{C3}), (\text{C4}) \text{ and } (\text{C5}) \text{ satisfied}\}]}{\mathsf{E}[\#\{t \in [0, 1]; W(vt, t) = 0, (\partial/\partial t) W(vt, t) > 0, W_x(vt, t) < 0\}]}.$$

The following theorem is the encountered counterpart of Theorem 7.1.

THEOREM 7.3. *Let $W : \mathbb{R}^2 \to \mathbb{R}$ be a stationary Gaussian process having a.s. twice continuously differentiable sample paths. Assume that the spectrum of $W$ has a continuous component and that the variance of the number of zeros $W_x(y, 0) = 0$, $y \in [0, 1]$, is finite. From $W$ define a process $\eta$ by the relation (7.8) and a set $B$ as in* (C1). *Furthermore, define a vector valued process $Y = (W, W_x, W_{xx})$, a set $\Gamma_{tx}$ as in* (C3) *and let $g(Y, x)$ be an indicator function defined by*

$$g(Y, x) = \mathbf{1}_{\{W_x(s,t) \le 0, \forall s \in \Gamma_{tx}\}} \mathbf{1}_{\{W(vt - x_2, t) \le u\}} \mathbf{1}_{\{W(vt + x_3, t) \le w\}}$$
$$\times \mathbf{1}_{\{W_{xx}(vt - x_2, t) < 0\}} \mathbf{1}_{\{W_{xx}(vt + x_3, t) > 0\}} \mathbf{1}_{\{(\partial/\partial t) W(vt, t) > 0\}}.$$

*Then, writing $N_0^\eta(B, g) = \#\{x \in B; \eta(x) = 0, g(Y, x) = 1\}$,*

$$\mathsf{E}[N_0^\eta(B, g)] = \int_B \mathsf{E}[|\det(\eta_x(x))| g(Y, x) \mid \eta(x) = 0] f_{\eta(x)}(0) \, dx,$$

*where both members are finite.*

The proof of the theorem follows the same lines as the proof of Theorem 7.1 and is, therefore, omitted.

We are now ready to state the main theorem for the distribution of $(x_2, x_3, H_2, H_3)$ in the encountered case. The proof is analogous to that of Theorem 7.2.

**Theorem 7.4.** *Let the assumptions and notation of Theorem 7.3 hold. Define*

(7.10) $\quad \beta(x) = Z_t(t)^+ W_{xx}(vt - x_2, t)^- W_{xx}(vt + x_3, t)^+ \mathbf{1}_{\{W_x(s,t) \leq 0, \forall s \in \Gamma_{tx}\}},$

*where $Z(t) = W(vt, t)$ is the encountered process.*

*Moreover, define a function $h_e(x,y) = h_e(x_1, x_2, x_3, y_1, y_2)$ by*

$$h_e(x,y) = \mathsf{E}[\beta(x) \mid \eta(x) = 0, W(vt - x_2, t) = y_1, W(vt + x_3, t) = y_2]$$
$$\times f_{\eta(x), W(vt-x_2,t), W(vt+x_3,t)}(0, y_1, y_2).$$

*Then the Palm distribution $F_{x_2, x_3, H_2, H_3}^{\mathrm{enc}}$ is given by*

(7.11)
$$F_{x_2, x_3, H_2, H_3}^{\mathrm{enc}}(r, s, u, w)$$
$$= \frac{\int_0^r \int_0^s \int_0^u \int_{-\infty}^w h_e(0, x_2, x_3, y_1, y_2) \, dy_2 \, dy_1 \, dx_3 \, dx_2}{\mathsf{E}[Z_t(0)^+ \mathbf{1}_{\{W_x(0,0)<0\}} \mid Z(0) = 0] f_{Z(0)}(0)},$$

*having the density*

(7.12)
$$f_{x_2, x_3, H_2, H_3}^{\mathrm{enc}}(r, s, u, w)$$
$$= \frac{h_e(0, r, s, u, w)}{\mathsf{E}[Z_t(0)^+ \mathbf{1}_{\{W_x(0,0)<0\}} \mid Z(0) = 0] f_{Z(0)}(0)},$$

*where $r, s, u > 0$ and $w < 0$.*

Earlier it was shown that the relation between spatial and encountered distribution of slope could be given a physical interpretation in terms of wave velocities. This is also true when it comes to the distribution of distances and heights of the local extrema closest to a wave center. To see that this is the case, consider the functions $\alpha(x_1, x_2, x_3)$ and $\beta(t, x_2, x_3)$ from Theorems 7.2 and 7.4, respectively. Due to a relation given in the proof of Corollary 6.1, namely $Z_t(0)^+ \mathbf{1}_{\{W_x(0,0)<0\}} = W_x(0,0)^- (V - v)^+$, the functions $\alpha$ and $\beta$ are related by

$$\beta(0, x_2, x_3) = (V - v)^+ \alpha(0, x_2, x_3),$$

where $V = -W_t(0,0)/W_x(0,0)$ is the wave velocity at the point $(x,t) = (0,0)$. Moreover, it is easy to see that $\eta(0, x_2, x_3) = \xi(0, x_2, x_3)$ so that, in fact, the functions $h_e(0, x_2, x_3, y_1, y_2)$ and $h_s(0, x_2, x_3, y_1, y_2)$ only differ by the velocity factor $(V - v)^+$ that enters into the expectation in the expression for $h_s$.

7.3. *Evaluation of the densities.* So far it has not been mentioned how the densities (7.7) and (7.12) can be evaluated. It turns out that for both densities the denominators can be computed explicitly. In the first, spatial, case, this is straightforward and

$$\mathsf{E}[W_x(0)^- \mid W(0) = 0] f_{W(0)}(0) = \frac{1}{2\pi} \sqrt{\frac{\lambda_{20}}{\lambda_{00}}}.$$

For the denominator of the encountered density (7.12) this is not as easy. However, [16] showed that by using certain symmetry properties of the Gaussian distribution one has

(7.13)
$$\mathsf{E}[Z_t(0)^+ \mathbf{1}_{\{W_x(0,0)<0\}} \mid Z(0) = 0] f_{Z(0)}(0)$$
$$= \frac{1}{4\pi} \sqrt{\frac{\lambda_{20}}{\lambda_{00}}} \left( \sqrt{v^2 + 2v \frac{\lambda_{11}}{\lambda_{20}} + \frac{\lambda_{02}}{\lambda_{20}}} - \frac{\lambda_{11}}{\lambda_{20}} - v \right),$$

whenever the spectral moments $\lambda_{00}$, $\lambda_{20}$, $\lambda_{02}$ and $\lambda_{11}$ are available.

For the numerators, on the other hand, there are no explicit formulas that can be used. Fortunately, there are very efficient numerical routines. In particular routines from the Matlab toolbox WAFO (Wave Analysis for Fatigue and Oceanography, available gratis at the web page http://www.maths.lth.se/matstat/wafo/), custom made for this type of calculations can be used; see [6]. However, before any computations can be done, the infinite dimensional indicators in the expressions for the densities must be approximated by an indicator on a finite set of grid points. In both cases, the indicator to be approximated is given by $\mathbf{1}_{\{W_x(u,0)\leq 0,-r\leq u\leq s\}}$; see (7.5) and (7.10). Therefore, let $\mathbf{U} = (u_1, \ldots, u_n)$, where $-r < u_1 < \cdots < u_n < s$, be a subdivision of the interval $[-r,s]$ and approximate the indicator $\mathbf{1}_{\{W_x(u,0)\leq 0,-r\leq u\leq s\}}$ by $\mathbf{1}_{\{W_x(u_i,0)\leq 0,i=1,\ldots,n\}}$. In this way, due to the inequality $\mathbf{1}_{\{W_x(u,0)\leq 0,-r\leq u\leq s\}} \leq \mathbf{1}_{\{W_x(u_i,0)\leq 0,i=1,\ldots,n\}}$, an upper bound of the densities is obtained. More information on computational details are given in the Appendix.

**8. Examples.** The methods presented in this paper can be applied to a Gaussian sea having a general spectrum. In the following examples, we will evaluate the distributions for the different wave characteristics for a longcrested Gaussian sea, that is, a Gaussian sea where the spectrum has no angular dependency. To get a longcrested sea let $S(\omega, \theta) = S(\omega)$, for one single fixed direction $\theta$. In this example $\theta = \pi$, meaning that the waves are moving in the direction of the positive $x$-axis. The spectrum $S(\omega)$ will be a JONSWAP frequency spectrum; see Figure 2. This is a family of spectra fully characterized by the set of parameters $(h_s, t_p, \gamma, \sigma_a, \sigma_b)$, where $h_s$ is the significant wave height, defined as four times the standard deviation of the sea elevation, and $t_p$ the peak period. The parameter $\gamma$, sometimes called the peak enhancement factor, determines the concentration of the spectrum around the peak frequency and $\sigma_a$ and $\sigma_b$ are spectral width parameters. Here, for illustration reasons only, we choose $h_s = 11.5$ m, $t_p = 12.25$ s, $\gamma = 1$, $\sigma_a = 0.07$ and $\sigma_b = 0.09$. Because $\gamma$ equals one this is also known as a Pierson–Moskowitz or Bretschneider spectrum. In this example, we have also specified a cut-off frequency $\omega_c$, that is, a frequency such that $S(\omega) = 0$ for all $\omega > \omega_c$, namely $\omega_c = 1.25$ rad/s.

8.1. *Distribution of slope.* The distribution of slope observed at centers of waves and overtaking encountered waves are given in Theorems 6.1 and 6.2, respectively. In Figure 3, these distributions are shown, in the encountered case for three different ship velocities. Clearly, the encountered distribution is shifted toward less steep waves and the faster the ship sails, the less steep the waves are. A simple explanation to this phenomenon is given by the dispersion relation for deterministic waves. Recall that by Corollary 6.1, the spatial distribution (6.3) is transformed into its encountered version (6.6) by the factor $(V - v)^+$. The physical interpretation of this factor is that waves that are likely to be slower than the ship will cancel. Because the dispersion relation (2.1) implies that steep deterministic waves are slow, it thus holds that mainly steep waves are canceled, leading to the shift in the distribution—the steep waves are simply not fast enough. However, to

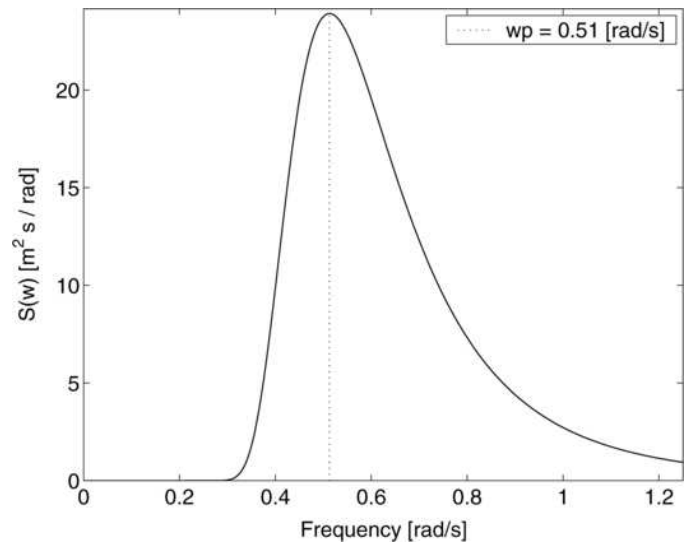

FIG. 2. *The JONSWAP frequency spectrum used in the examples. This spectrum has significant wave height $h_s = 11.5$ m, peak period $t_p = 12.25$ s and cut-off frequency $\omega_c = 1.25$ rad/s.*

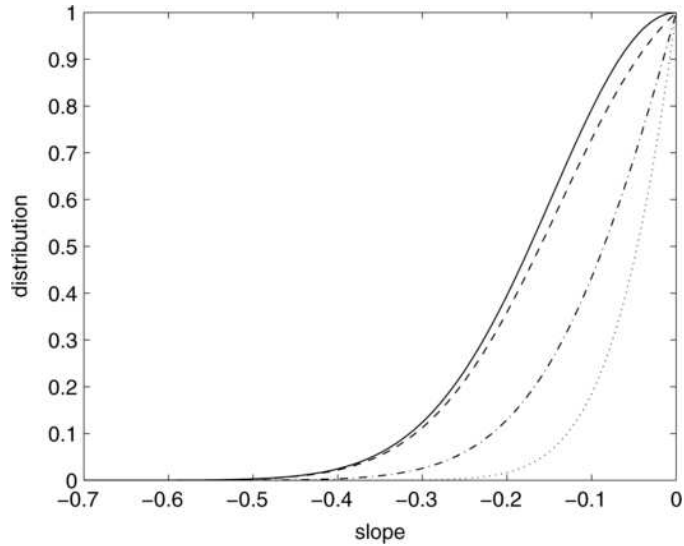

Fɪɢ. 3.   *Palm distribution of slope observed at centers of waves in space (solid) and at encountered centers of waves overtaking a ship sailing with velocity 7 m/s (dashed),* 13 m/s *(dash-dotted) and* 16 m/s *(dotted). The computations are done for a JONSWAP frequency spectrum having significant wave height* $h_s = 11.5$ *m, peak period* $t_p = 12.25$ *s and cut-off frequency* $\omega_c = 1.25$ *rad/s.*

fully understand the relation between slope and velocity, it is too great a simplification to look at deterministic waves, although it gives some insight into the physics involved in the problem. For a more thorough understanding, one should study the distribution of the random velocity of the waves conditional on the slope.

8.2. *Distribution of waveheight and half-wavelength.*   Next we consider the spatial and encountered joint densities of the waveheight $H_2 - H_3$ and the half-wavelength $x_2 + x_3$, that can be obtained from Theorems 7.2 and 7.4. In Figure 4 these densities are shown, in the encountered case for two different ship velocities, namely $v = 7, 13$ m/s. Clearly, the effect of observing the sea surface at waves overtaking the vessel compared to observing it along a line at a fixed time, is that the density is shifted toward longer waves in the former case. Moreover, the shift increases with increasing ship velocity. As for the distribution of slope, this behavior can to some extent be explained by the dispersion relation. Recall that [see (4.1)] the velocity of deterministic waves satisfying the dispersion relation is proportional to the square root of the wavelength. This means that long waves are faster than short ones. Because the factor $(V - v)^+$, explaining the difference between the spatial and encountered densities, eliminates waves that are slower than the ship,

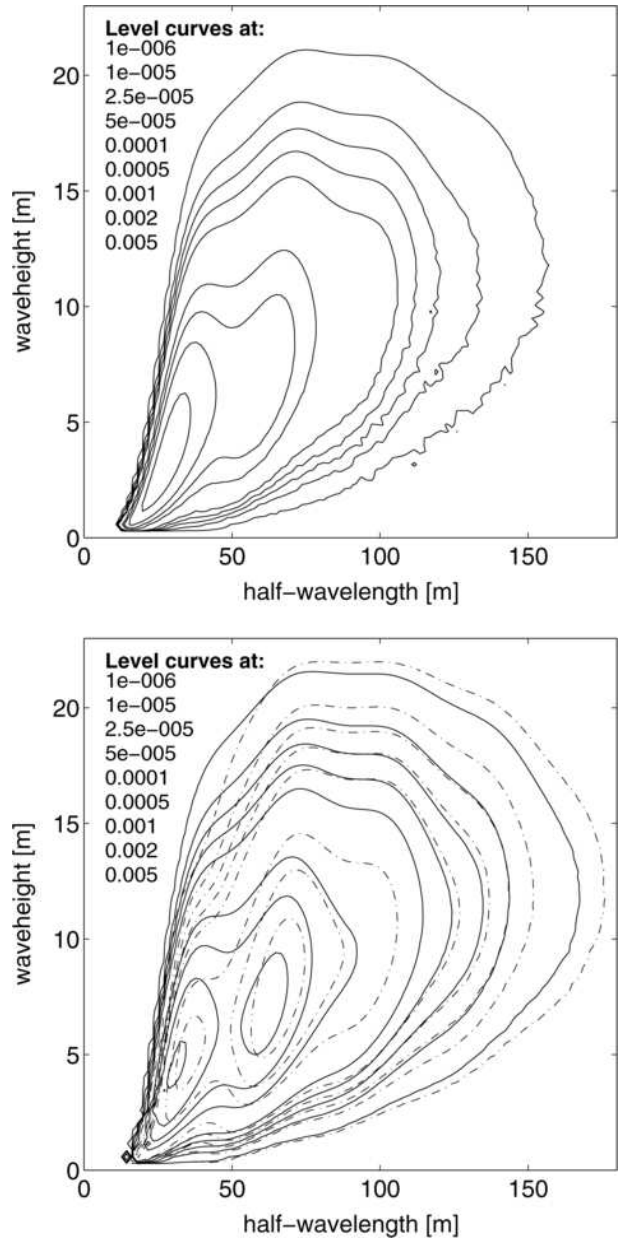

FIG. 4. *Density of half-wavelength $x_2 + x_3$ and waveheight $H_2 - H_3$ observed at wave centers along a line at a fixed time point (top) and at centers of waves overtaking the ship (bottom). In the latter case, the velocities of the ship are $v = 7$ m/s (solid) and $v = 13$ m/s (dash-dotted). The computations are done for a JONSWAP frequency spectrum having significant wave height $h_s = 11.5$ m, peak period $t_p = 12.25$ s and cut-off frequency $\omega_c = 1.25$ rad/s.*

it thereby eliminates small and moderately high waves in favor of long and high waves.

## APPENDIX: COMPUTATIONAL DETAILS

**A.1. Evaluation of densities with WAFO.** In this section we discuss the evaluation of the densities $f^{\text{space}}$ and $f^{\text{enc}}$ from Theorems 7.2 and 7.4, respectively, using the WAFO-toolbox.

The density $f^{\text{space}}$ in (7.7) can be computed by the routine `spec2mmtpdf` and it computes the desired density given the spectrum of the process $W(x, t)$. For the encountered density $f^{\text{enc}}$ given by (7.12), however, there is no existing routine in the toolbox, but the function `rind` that evaluates Gaussian multivariate expectations is of great help.

Before the intensity in the numerator in (7.12) can be evaluated, the infinite dimensional indicator must be approximated by an indicator on a finite set of grid points of the interval $[-r, s]$. Let $\mathbf{U} = (u_1, \ldots, u_n)$, where $-r < u_1 < \cdots < u_n < s$, be a subdivision of the interval $[-r, s]$ and approximate the indicator $\mathbf{1}_{\{W_x(u, 0) < 0, -r \leq u \leq s\}}$ by $\mathbf{1}_{\{W_x(u_i, 0) < 0, i = 1, \ldots, n\}}$. To use `rind` the full distribution, that is, the means and covariances of all variables involved, must be computed. To simplify notation the variables are collected in the following vectors:

$$X_t = (W_x(u_1, 0), \ldots, W_x(u_n, 0)) = W_x(\mathbf{U}, 0),$$

$$X_d = (W_{xx}(-r, 0), W_{xx}(s, 0), vW_x(0, 0) + W_t(0, 0)),$$

$$X_c = (W_x(-r, 0), W_x(s, 0), W(-r, 0), W(s, 0), W(0, 0)).$$

Because $W(x, t)$ is a zero mean stationary process, and due to the fact that the means of the (quadratic mean) derivatives of a stationary process are zero, the mean of all three vectors above is zero.

Computing the covariances requires more effort. In general, for a stationary field $W(\mathbf{t})$ with covariance function $R(\tau)$ the following differentiation rule holds

$$(1.14) \quad \text{Cov}\left(\frac{\partial^{a+b} W(\mathbf{s})}{\partial^a s_j \partial^b s_k}, \frac{\partial^{c+d} W(\mathbf{t})}{\partial^c t_m \partial^d t_n}\right) = (-1)^{a+b} \frac{\partial^{a+b+c+d} R(\tau)}{\partial^a \tau_j \partial^b \tau_k \partial^c \tau_m \partial^d \tau_n}\bigg|_{\tau = \mathbf{t} - \mathbf{s}}.$$

Let $\widetilde{\Sigma}$ be the covariance matrix of $(X_t, X_d, X_c)$ and partition it as

$$\widetilde{\Sigma} = \begin{pmatrix} \Sigma_{tt} & \Sigma_{td} & \Sigma_{tc} \\ & \Sigma_{dd} & \Sigma_{dc} \\ & & \Sigma_{cc} \end{pmatrix},$$

where it is understood that the lower left part of the matrix is the transpose of the upper right part. By repeated use of (1.14) one obtains

$$\Sigma_{tt} = (-R_{\xi\xi}(\mathbf{U}^T - \mathbf{U}, 0)),$$

$$\Sigma_{td} = (-R_{\xi\xi\xi}(-r - \mathbf{U}^T, 0), -R_{\xi\xi\xi}(s - \mathbf{U}^T, 0),$$
$$- vR_{\xi\xi}(-\mathbf{U}^T, 0) - R_{\xi\tau}(-\mathbf{U}^T, 0)),$$
$$\Sigma_{tc} = (-R_{\xi\xi\xi}(-r - \mathbf{U}^T, 0), -R_{\xi\xi}(s - \mathbf{U}^T, 0), -R_{\xi}(-r - \mathbf{U}^T, 0), \ldots,$$
$$- R_{\xi}(s - \mathbf{U}^T, 0), -R_{\xi}(-\mathbf{U}^T, 0)),$$

where $-R_{\xi\xi}(\mathbf{U}^T - \mathbf{U}, 0)$ should be interpreted as the matrix with element on row $i$ and column $j$ equal to $-R_{\xi\xi}(u_i - u_j, 0)$, and $-R_{\xi\xi\xi}(-r - \mathbf{U}^T, 0)$ is a column vector with elements $-R_{\xi\xi\xi}(-r - u_j, 0)$, $i, j = 1, \ldots, n$. The other submatrices should be interpreted in a similar fashion. The remaining part of $\Sigma$ is equal to

$$\Sigma_{dd} = \begin{pmatrix} R_{\xi\xi\xi\xi}(0,0) & R_{\xi\xi\xi\xi}(s+r,0) & vR_{\xi\xi\xi}(r,0) + R_{\xi\xi\tau}(r,0) \\ & R_{\xi\xi\xi\xi}(0,0) & vR_{\xi\xi\xi}(-s,0) + R_{\xi\xi\tau}(-s,0) \\ & & -v^2 R_{\xi\xi}(0,0) - 2vR_{\xi\tau}(0,0) - R_{\tau\tau}(0,0) \end{pmatrix},$$

$$\Sigma_{dc} = \begin{pmatrix} 0 & R_{\xi\xi\xi}(s+r,0) & \\ R_{\xi\xi\xi}(-r-s,0) & 0 & \cdots \\ -vR_{\xi\xi}(-r,0) - R_{\xi\tau}(-r,0) & -vR_{\xi\xi}(s,0) - R_{\xi\tau}(s,0) & \end{pmatrix}$$

$$\cdots \begin{pmatrix} R_{\xi\xi}(0,0) & R_{\xi\xi}(s+r,0) & R_{\xi\xi}(r,0) \\ R_{\xi\xi}(-r-s,0) & R_{\xi\xi}(0,0) & R_{\xi\xi}(-s,0) \\ -vR_{\xi}(-r,0) - R_{\tau}(-r,0) & -vR_{\xi}(s,0) - R_{\tau}(s,0) & 0 \end{pmatrix},$$

$$\Sigma_{cc} = \begin{pmatrix} -R_{\xi\xi}(0,0) & -R_{\xi\xi}(s+r,0) & 0 & -R_{\xi}(s+r,0) & -R_{\xi}(r,0) \\ & -R_{\xi\xi}(0,0) & -R_{\xi}(-r-s,0) & 0 & -R_{\xi}(-s,0) \\ & & R(0,0) & R(s+r,0) & R(r,0) \\ & & & R(0,0) & R(-s,0) \\ & & & & R(0,0) \end{pmatrix}.$$

The covariance function and its derivatives can be computed by the WAFO-function `spec2cov`, that from a given spectrum computes the covariance function and its derivatives up to the fourth order.

**A.2. Computational issues.** When the conditional expectation in the numerator of (7.12) is computed numerically new difficulties arise. If the covariance matrix is badly scaled, densities can become singular in the numerical computations. For example, it might happen that variances become negative even though they should be small and positive. A way to make the algorithms more stable is to rescale the covariance matrix such that its elements are of comparable sizes.

For the spectra we will consider the average period is about 10 seconds and the value of the average wave length is about 10 times as large, but measured in meters. Because the sea elevation takes values in the same range in both time and space, the derivatives will be of very different magnitudes.

To get a better scaled covariance matrix the time- and space-coordinates are, therefore, transformed according to

$$\tilde{x} = x\sqrt{\frac{\lambda_{20}}{\lambda_{00}}}, \qquad \tilde{t} = t\sqrt{\frac{\lambda_{02}}{\lambda_{00}}}.$$

Making this change of variables the spectral density is transformed to

$$\tilde{S}(\tilde{\omega}, \theta) = S(\omega, \theta)\frac{\sqrt{\lambda_{02}}}{\lambda_{00}^{3/2}}.$$

Note that the angle $\theta$ remains unaffected. In this way, the spectral moments $\tilde{\lambda}_{00}$, $\tilde{\lambda}_{20}$ and $\tilde{\lambda}_{02}$ of $\tilde{S}(\tilde{\omega}, \theta)$ all become equal to one; see [17] for a proof. This means that the average number of waves per transformed meter is equal to $1/2\pi$ as well as the number of waves per transformed second. Also, the variance of the sea elevation is transformed to one. Note that it is common practice in wave data analysis to normalize the data to have unit variance. However, in this case this is not enough in order to avoid a badly scaled covariance matrix because the problems arise due to different spatial and temporal scales. When the coordinates are transformed, the velocity of the ship changes to

$$\tilde{v} = v\sqrt{\frac{\lambda_{20}}{\lambda_{02}}}.$$

To relate the spectral moments to the elements in the covariance matrix note that according to (2.2) and (2.3) it holds that $\lambda_{00} = R(0,0)$, $\lambda_{20} = -R_{\xi\xi}(0,0)$ and $\lambda_{02} = -R_{\tau\tau}(0,0)$. Because many of the elements on the diagonal of $\widetilde{\Sigma}$ are of this form, we can conclude that the covariance matrix obtained after the coordinate transformation is better scaled than the original one. After the density is computed with the transformed coordinates, it is then easily transformed back to the true units seconds and meters. These kind of transformations are implemented in the function wnormspec in WAFO.

S. ABERG
I. RYCHLIK
MATHEMATICAL SCIENCES
CHALMERS UNIVERSITY OF TECHNOLOGY
AND
GÖTEBORG UNIVERSITY
SE-412 96 GÖTEBORG
SWEDEN
E-MAIL: abergs@chalmers.se
       rychlik@chalmers.se

R. LEADBETTER
DEPARTMENT OF STATISTICS
  AND OPERATIONS RESEARCH
UNIVERSITY OF NORTH CAROLINA
CHAPEL HILL, NORTH CAROLINA 27599
USA
E-MAIL: mrl@email.unc.edu